\documentclass[12pt]{article}
\usepackage{amsmath}                 
\usepackage{amssymb}                 
\usepackage{amsthm}                  

\newcommand{\GKdim}{\mathop{\rm GKdim}\nolimits}
\newcommand{\Bfield}{\mathfrak K}

\newtheorem{definition}{Definition}
\newtheorem{lemma}{Lemma}
\newtheorem{proposition}{Proposition}

\begin{document}

\title{ Simple algebra with arbitrary\\
        odd Gel'fand-Kirillov dimension }
\author{S.~S. Konyuhov}
\date{}
\maketitle

\begin{abstract}
It is given an example of finitely generated (f.~g.) simple algebra 
over a field $\Bfield$ (char $\Bfield=0$) with arbitrary odd Gel'fand-Kirillov dimension. 
\end{abstract}

\section{Introduction}

Let $\mathcal A -$ algebra over a field $\Bfield,$ $V -$ its finite-dimensional subspace containing the 
unit $1_{\mathcal A}\in\mathcal A$ and generating $\mathcal A$ over $\Bfield.$ Then there is an increasing 
chain of subspaces:
$$
  \Bfield \subseteq V \subseteq V^2 \subseteq\cdots\subseteq V^{n}
    \subseteq\cdots\subseteq\bigcup\limits_{n=0}^{\infty} V^{n} =\mathcal A,\,
     \dim_{\Bfield}(V^{n}) < \infty, \,\forall n\in\mathbb N_{0}.
$$

Asymptotic behaviour of this sequence is the invariant known as \emph{growth} or 
\emph{Gel'fand-Kirillov dimension} and defined by formula:
$$
\GKdim(\mathcal A) = \limsup\limits_{n\to\infty}\frac{\log\dim_{\Bfield}(V^{n})}{\log n}. 
$$

It was introduced by I.~M. Gel'fand and A.~A. Kirillov in 1966~\cite{paper_GK1, paper_GK2}. In these works 
they have made the conjecture known now as Gel'fand-Kirillov hypothesis and concerning problems of Lie 
algebra isomorphisms.

In~\cite{paper_Milnor} Milnor has given the corresponding notion for groups and showed that the growth 
of the fundumental group of Riemann manifold is connected with its curvature. More recent reference 
to this connection is mentioned in Schwartz' work~\cite{paper_Schwartz}.

First systematical study of the Gel'fand-Kirillov dimension was made by Borho and Kraft in 
1976~\cite{paper_BK}. However, this work contains some errors found later. Since  1976 it 
appears a lots of articles utilizing properties of the Gel'fand-Kirillov dimension. Now it is 
one of the standard invariants for study of noncommutative algebras.

Previously for these purposes the Krull dimension was used but it turns out the most of theorems need
additional technical hypotheses such as symmetry and invariance of ideals. Gel'fand-Kirillov dimension 
satisfies these two conditions and in the most cases it may be calculated easier. Sometimes knowledge of 
the Gel'fand-Kirillov dimension gives information about the Krull dimension.

The Gel'fand-Kirillov dimension of finite-dimensional algebra equals 0 and so it measures deviation of 
 arbitrary algebra $\mathcal A$ from the finite-dimensional case.

Commutativity and algebraic dependency of algebra generators have special influence on the 
Gel'fand-Kirillov dimension.

For example, if $\mathcal A -$ f.~g. commutative integral domain then  $\GKdim(\mathcal A)$ coincides 
with the transcendence degree of its fraction field and so equals to the number of independent variables
in the largest polynomial algebra included in $\mathcal A.$ 

And conversely, a free algebra with two generators has infinite Gel'fand-Kirillov dimension.
 
As one of the most general problems dealing with Gel'fand-Kirillov dimension it is often considered
the question whether f.~g. algebra from given class and having arbitrary Gel'fand-Kirillov
dimension exists. For some algebra classes this problem was solved positevily.
 
The Gel'fand-Kirillov dimension of f.~g. commutative algebra is integer~\cite{paper_Dr}.
Relatively free algebras have integer Gel'fand-Kirillov dimension~\cite{paper_Dr} too. In \cite{paper_I}
an example of semiprime algebra with arbitrary integer Gel'fand-Kirillov dimension was build.

In \cite{paper_IW} Irving and Warfield gave an example of primitive algebra with arbitrary 
Gel'fand-Kirillov dimension assuming existence of finite-dimensional algebraic extension of basic field.
Author of~\cite{paper_UV} has constructed an example of primitive algebra with arbitrary Gel'fand-Kirillov 
dimension over arbitrary field. 

\newpage
For PI-algebras it is known~\cite{paper_BK, paper_KL}:
\begin{enumerate}
\item if $\mathcal A -$ prime PI-algebra then $\GKdim(\mathcal A)\in\mathbb N;$ 
\item if $\mathcal A -$ f.~g. simple PI-algebra then $\GKdim(\mathcal A)=0;$
\item if $\mathcal A -$ noetherian PI-algebra then $\GKdim(\mathcal A)\in\mathbb N;$
\item $\forall\gamma\in\mathbb R, \gamma\geqslant 2~\exists
      \mbox{ f. g. PI-algebra with } \GKdim(\mathcal A)=\gamma.$
\end{enumerate}

Weyl algebra $\mathcal A_{w}$ is the example of simple algebra with arbitrary even Gel'fand-Kirillov 
dimension $\GKdim(\mathcal A_{w})=2n, n\in\mathbb N.$ In \cite{paper_SSW} it was shown that algebra with 
$\GKdim(\mathcal A)=1$ is PI-algebra therefore there is no simple algebra with $\GKdim(\mathcal A)=1.$ 
In this work an example of simple algebra with $\GKdim(\mathcal A_{w})=2n+1, n\in \mathbb N$ is given.

\newpage

\section{ Finitely generated simple algebra with \\
          arbitrary odd Gel'fand-Kirillov dimension }

\begin{definition} 
Let $\Bfield$ be a field with $char(\Bfield) = 0$ then define $\Bfield$-algebra $\mathcal A$ as follows:
$\mathcal A = \Bfield\left\langle x, y, z \mid [x, y] = [y, z] = 1, [x, z] = -xy-yz \right\rangle
\footnote{As usual $[x, y] = xy - yx$ and so on.}.$
\end{definition}

First we check that our definition is correct and $\mathcal A$ is nonempty associative $\Bfield$-algebra. 
Indeed, it can be considered as factor-algebra $\Bfield\langle x, y, z\rangle/I$ where 
$I =(f_{1}, f_{2}, f_{3}), f_{1} = yx - (xy-1) = \omega_{1}- g_{1}, f_{2} = zy-(yz-1) = \omega_{2}- g_{2},\\ f_{3} = zx - (xz+xy+yz)=\omega_{3}- g_{3}.$

If we introduce the deg-lex ordering on $\mathcal A$ by relations: $x \prec y \prec z$ then polynoms 
$f_{1}, f_{2}$ form the only critical pair because of $z\omega_{1} = z(yx) = (zy)x = \omega_{2}x$ 
and the according critical $s$-polynom is $s_{2, 1} = f_{2}x - zf_{1}.$ 

Have defined reduction operators for a polynom system $\mathcal F = \{f_{1}, f_{2}, f_{3}\}:$
$$
r_{u i v}(u\omega_{i}v) \triangleq ug_{i}v, \;
r_{u i v}(w) \triangleq  w, \; w \ne u\omega_{i}v, \; u, v, w \in \langle x, y, z \rangle ,
$$
we get that $s_{2,1}$ is reduced to 0:
\begin{multline*}
s_{2, 1} = f_{2}x-zf_{1} = \left((zy-(yz-1)\right)x - z\left(yx-(xy-1)\right)\\
= (zy)x - (yz)x + x- z(yx) + z(xy) - z =  -y(xz+xy+yz) + x + (xz+xy+yz)y - z\\
= -(xy-1)z - (xy-1)y - y^{2}z + x + x(yz-1) + xy^{2} + y(yz-1) - z\\
= -xyz + z - xy^{2} + y - y^{2}z + x + xyz - x + xy^{2} + y^{2}z -y - z = 0.
\end{multline*}
So $\mathcal F$ is a Gr\"{o}bner basis of the polynomial ideal $I.$ And so each element of $\mathcal A$ 
has unique normal form i.~e.:
$
\forall w \in \mathcal A \quad\exists! \alpha_{ijk}\in \Bfield : w = \sum\alpha_{ijk} x^{i}y^{j}z^{k}.
$
Moreover by definition we have $\{x^{i}y^{j}z^{k} \mid i, j, k\in\mathbb N_{0}\}-$ a basis of 
$\mathcal A$ over $\Bfield.$ 

Now let $V_{0} = \Bfield$ and $V = V_{1} = \Bfield x + \Bfield y + \Bfield z$ 
then the set $\mathcal S=\{x^{i}y^{j}z^{k} \mid 0 \leqslant i+j+k\leqslant n\}$ is a basis of linear 
subspace $V_{n} = V + V^{2} + \ldots + V^{n}$ forming filtration of 
$\mathcal A=\bigcup\limits_{n=0}^{\infty}V_{n}.$

So bearing in mind that: 
$
\dim_{\Bfield} V_{n} = \sum\limits_{i=0}^{n} \frac{(i+1)(i+2)}{2} =
        \frac{1}{2}\sum\limits_{i=0}^n (i^2+3i+2) =
        \frac{n^3}6+n^2+\frac{11}{6}n,
$
we yield:
$$
\GKdim(\mathcal A) = 
   \limsup\limits_{n\to\infty}\frac{\log\dim_{\Bfield}(V_{n})}{\log n} =
   \lim\limits_{n\to\infty}\frac{3\log n+\log(\frac{1}{6}+\frac{1}{n}+\frac{11}{6n^2})}{\log n} = 3.
$$
As a result we have proved
\begin{proposition}\label{prop_dimension}
$\GKdim(\mathcal A) = 3.$
\end{proposition}

The next lemma will be used throughout all the rest of this report.
\begin{lemma}\label{lemma_AlgebraRelations}
For elements of the algebra $\mathcal A$ the following relations hold
\footnote{We assume always that $kx^{k-1}, ky^{k-1}$ and $kz^{k-1}$ equal 0 if $k=0.$}:
\begin{align*}
& [x^{k}, y] = kx^{k-1},\quad [x, y^{k}] = ky^{k-1},\quad [y, z^{k}] = kz^{k-1},\quad [y^{k}, z] = ky^{k-1}, \\
   & [x, (x+z)^{k}] = -k\left(y(x+z)^{k} + (x+z)^{k-1}\right) = [(x+z)^{k}, z],                  \\
   & [(y+z)^{k}, z] = k(y+z)^{k-1},\quad [y,(x+z)^k] = 0, \quad \forall k \in\mathbb N_{0}.
\end{align*}
\end{lemma}
\begin{proof}
Induction by $k.$

Prove e.~g. that $[x, (x+z)^{k}] = -k\left(y(x+z)^{k} + (x+z)^{k-1}\right).$
By $k=1$ we have: $[x, x+z]=[x,x]+[x,z]=-yz-(yx+1)=-y(x+z)-1.$ 

If now $k>1$ then:
\begin{multline*}
  (x+z)^{k+1}x = (x+z)(x+z)^{k}x = \\
 = (x+z)\left(x(x+z)^{k} + ky(x+z)^{k} + k(x+z)^{k-1}\right) = \\
 = \left(x(x+z) + y(x+z) + 1\right)(x+z)^{k} + ky((x+z)^{k+1} + k(x+z)^{k} = \\
 = k(x+z)^{k+1} + (k+1)y(x+z)^{k+1}+(k+1)(x+z)^{k}.
\end{multline*}

So $[x, (x+z)^{k+1}] = -(k+1)\left(y(x+z)^{k+1}+(x+z)^k\right).$

\end{proof}

For future we need the following
\begin{lemma}
$\forall N, M, K$ and $L \in\mathbb N_{0}:$
$$
 y^{N}(x+z)^{K} = y^{M}(x+z)^{L} \Longrightarrow N=M \mbox{ and } K = L.
$$
\end{lemma}
\begin{proof}
First show by the induction that $\forall k\in\mathbb N:$ 
$$
(x+z)^{k} = \sum\limits_{i+j=k}\alpha_{ij}x^{i}z^{j} + 
            \sum\limits_{p+q+r\leqslant k}\beta_{pqr}x^{p}y^{q}z^{r}
$$ 
where $\alpha_{ij}>0$ and $\beta_{p 0 k-p}=0.$

This is clear if $k=1.$ So let $k>1$ then using lemmas~\ref{lemma_AlgebraRelations} and~\ref{lemma_zKx} 
we obtain:
\begin{multline*}
(x+z)^{k+1} = (x+z)^{k}(x+z) =\\
= \sum\limits_{i+j=k}\alpha_{ij}x^{i}z^{j}x + 
  \sum\limits_{p+q+r\leqslant k}\beta_{pqr}x^{p}y^{q}z^{r}x + 
  \sum\limits_{i+j=k}\alpha_{ij}x^{i}z^{j+1} + 
  \sum\limits_{p+q+r\leqslant k}\beta_{pqr}x^{p}y^{q}z^{r+1} =\\
= \sum\limits_{i+j=k}\alpha_{ij}x^{i}\left(\sum\limits_{s+t=j}\bar{\alpha}_{st}xy^{s}z^{t} + 
  \sum\limits_{s+t<j-1}\bar{\beta}_{st}xy^{s}z^{t} + 
  \sum\limits_{s+t\leqslant j+1}\bar{\gamma}_{st}y^{s}z^{t}\right) +\\
+ \sum\limits_{p+q+r\leqslant k}\beta_{pqr}x^{p}y^{q}\left(\sum\limits_{u+v=r}\tilde{\alpha}_{uv}xy^{u}z^{v}+
  \sum\limits_{u+v<r-1}\tilde{\beta}_{uv}xy^{u}z^{v} + 
  \sum\limits_{u+v\leqslant r+1}\tilde{\gamma}_{uv}y^{u}z^{v}\right)+ \\
+ \sum\limits_{i+j=k}\alpha_{ij}x^{i}z^{j+1} + 
  \sum\limits_{p+q+r\leqslant k}\beta_{pqr}x^{p}y^{q}z^{r+1} =\\
= \sum\limits_{i+j=k}\alpha_{ij}\bar{\alpha}_{0j}x^{i+1}z^{j} +
\sum\limits_{i+j=k}\sum\limits_{s+t=j,\atop s\ne 0}\alpha_{ij}\bar{\alpha}_{st}x^{i+1}y^{s}z^{t} + 
\sum\limits_{i+j=k}\sum\limits_{s+t<j-1}\alpha_{ij}\bar{\beta}_{st}x^{i+1}y^{s}z^{t} +\\
+ \sum\limits_{i+j=k}\sum\limits_{s+t\leqslant j+1}\alpha_{ij}\bar{\gamma}_{st}x^{i}y^{s}z^{t} + 
\sum\limits_{p+q+r\leqslant k}\sum\limits_{u+v=r}\beta_{pqr}\tilde{\alpha}_{uv}x^{p}(xy^{q}-qy^{q-1})y^{u}z^{v} +\\
+ \sum\limits_{p+q+r\leqslant k}\sum\limits_{u+v<r-1}\beta_{pqr}\tilde{\beta}_{uv}x^{p}(xy^{q}-qy^{q-1})y^{u}z^{v} +
\sum\limits_{p+q+r\leqslant k}\sum\limits_{u+v\leqslant r+1}\beta_{pqr}\tilde{\gamma}_{uv}x^{p}y^{q+u}z^{v} +\\
+ \sum\limits_{i+j=k}\alpha_{ij}x^{i}z^{j+1} +
\sum\limits_{p+q+r\leqslant k}\beta_{pqr}x^{p}y^{q}z^{r+1} = 
\sum\limits_{i+j=k+1}\hat{\alpha}_{ij}x^{i}z^{j} +
\sum\limits_{p+q+r\leqslant k+1}\hat{\beta}_{pqr}x^{p}y^{q}z^{r}.
\end{multline*}
Due to the induction we have $\beta_{p 0 k-p}=0$ and $\bar{\gamma}_{0 k+1}=0$ 
(lemma~\ref{lemma_zKx}, see later) therefore $\hat{\alpha}_{0 k+1} = \alpha_{0 k},$
$\hat{\alpha}_{k+1 0} = \alpha_{k 0},$ $\hat{\alpha}_{i j}=
\alpha_{i-1 j}\bar{\alpha}_{0 j} + \alpha_{i j-1}, i+j=k+1,$ $i\ne 0,k+1$ and
$\hat{\alpha}_{ij}>0, i+j=k+1$ as $\bar{\alpha}_{0 j}>0$ (lemma~\ref{lemma_zKx}) 
and $\alpha_{i j}>0, i+j=k$ (inductive supposition).

Further $\hat{\beta}_{p 0 k+1-p}=0$ for all monoms $x^{p}z^{k+1-p}$ with nonzero coefficients
were collected in the first sum.

Now using relations from lemma~\ref{lemma_AlgebraRelations} it is easy to write down the normal 
forms of $y^{N}(x+z)^{K}$ and $y^{M}(x+z)^{L}:$
\begin{gather*}
y^{N}(x+z)^{K} = \sum\limits_{i+j=K}\alpha_{ij}x^{i}y^{N}z^{j} +
                 \sum\limits_{p+q+r=K}\beta_{pqr}x^{p}y^{N+q}z^{r} +
                 \sum\limits_{p+q+r<K+N,\atop p+r<K}\gamma_{pqr}x^{p}y^{q}z^{r},\\
y^{M}(x+z)^{K} = \sum\limits_{i+j=L}\bar{\alpha}_{ij}x^iy^{M}z^{j} +
                 \sum\limits_{p+q+r=L}\bar{\beta}_{pqr}x^{p}y^{M+q}z^{r} +
                 \sum\limits_{p+q+r<L+M,\atop p+r< L}\bar{\gamma}_{pqr}x^{p}y^{q}z^{r}.
\end{gather*}
And one may check that different pairs $(N, K)$ and $(M,L)$ correspond to different polynoms 
$y^{N}(x+z)^{K}$ and $y^{M}(x+z)^{L}.$
\end{proof}

\textsc{Corollary.}
Polynoms $\{y^{i}(x+z)^{k} \mid i, k\in\mathbb N_{0}\}$ are linear independent over $\Bfield.$

In the proof given above we have used the supplementary
\begin{lemma}~\label{lemma_zKx}
$\forall k \in \mathbb N:$
$$
  z^{k}x = \sum\limits_{i+j=k}\alpha_{ij}xy^{i}z^{j} + \sum\limits_{p+q<k-1}\beta_{pq}xy^{p}z^{q} +
           \sum\limits_{i+j\leqslant k+1}\gamma_{ij}y^{i}z^{j}
$$
where $\alpha_{i j}>0$ for $i \ne 0, k+1,$ $\gamma_{0 k+1}=\gamma_{k+1 0}>0.$ 
\end{lemma}

\begin{proof} 
Divide the proof on two steps. First we show that $\forall k\in\mathbb N:$
$$
(y+z)^{k} = 
\sum\limits_{i+j=k}\bar{\alpha}_{ij}y^{i}z^{j} +
\sum\limits_{p+q<k-1}\bar{\beta}_{pq}y^{p}z^{q}, \mbox{ and } \bar{\alpha}_{ij}>0.
$$

This relation is clear if $k=1$ so let $k>1$ then: 
\begin{multline*}
 (y+z)^{k+1} = (y+z)^k(y+z) = \sum\limits_{i+j=k}\alpha_{ij}y^{i}(yz^{j}-jz^{j-1}) +\\
+ \sum\limits_{i+j=k}\alpha_{ij}y^{i}z^{j+1} + 
 \sum\limits_{p+q<k}\beta_{pq}y^{p}(yz^{q}-qz^{q-1}) + \sum\limits_{p+q<k}\beta_{pq}y^{p}z^{q-1} =\\
 = \sum\limits_{i+j=k+1}\bar{\alpha}_{ij}y^{i}z^{j} + \sum\limits_{p+q<k+1}\bar{\beta}_{pq}y^{p}z^{q} 
\end{multline*}
and for $i, j \ne 0, \ne k+1$ $\bar{\alpha}_{ij}>0$ since $\bar{\alpha}_{0 k+1}=\bar{\alpha}_{0 k},$ 
$\bar{\alpha}_{k+1 0}=\bar{\alpha}_{k 0},$ $\bar{\alpha}_{ij}=\alpha_{i-1 j}+\alpha_{i j-1}$ if 
$i, j \ne 0, k+1$ and $\alpha_{ij}>0$ by inductive supposition.

We finish the proof if we show that $\forall k\in\mathbb N:$
$$
z^{k}x = x(y+z)^{k} + \sum\limits_{i+j\leqslant k+1}\bar{\bar{\alpha}}_{ij}y^{i}z^{j}
$$ 
where $\alpha_{0 k+1}=\alpha_{k+1 0}=0,$ $\alpha_{i k+1-i}>0,$ $i\ne 0, k+1.$ 

Again we are using induction by $k.$ Case $k=1$ is evident so let $k \geqslant 1$ then:
\begin{multline*}
 z^{k+1}x = z\bigl(x(y+z)^{k} + \sum\limits_{i+j\leqslant k+1}\alpha_{ij}y^{i}z^{j}\bigr) = 
 (x(y+z)+yz)(y+z)^{k} +\\ 
+ \sum\limits_{i+j\leqslant k+1}\alpha_{ij}(y^{i}z - iy^{i-1})z^{j} = 
x(y+z)^{k+1} + y\bigl((y+z)^{k}z - k(y+z)^{k-1}\bigr) +\\
+ \sum\limits_{i+j\leqslant k+1}\alpha_{ij}y^{i}z^{j+1} -
\sum\limits_{i+j\leqslant k+1}\alpha_{ij}iy^{i-1}z^{j} = 
x(y+z)^{k+1} + \sum\limits_{i+j\leqslant k+1}\alpha_{ij}y^{i}z^{j+1} +\\
+ \sum\limits_{i+j=k}\bar{\alpha}_{ij}y^{i+1}z^{j+1} + 
\sum\limits_{p+q<k-1}\bar{\beta}_{pq}y^{p+1}z^{q+1}  - 
k\sum\limits_{i+j = k-1}{\tilde{\alpha}_{ij}}y^{i+1}z^{j} -\\
- k\sum\limits_{p+q<k-2}\tilde{\beta}_{pq}y^{p+1}z^{q} -
\sum\limits_{i+j\leqslant k+1}\alpha_{ij}iy^{i-1}z^{j} =\\
 = x(y+z)^{k+1} + \sum\limits_{i+j\leqslant k+2}\bar{\bar{\alpha}}_{ij}y^{i}z^{j}
\end{multline*}
where $\bar{\bar{\alpha}}_{0 k+2} = \alpha_{0 k+1} = 0,$ $\bar{\bar{\alpha}}_{k+2 0} = \alpha_{k+1 0} =0$
and $\bar{\bar{\alpha}}_{i j} = \bar{\alpha}_{i-1 j-1} + \alpha_{i j-1} >0,$ $i+j=k+2$ since
$\alpha_{i j-1} > 0$ by the induction and $\bar{\alpha}_{i-1 j-1}>0$ as we have obtained in the first step.
\end{proof}

\begin{definition}
Define linear operator $R_{y} \in End_{R}(\mathcal A)$ as: 
$$ R_{y}(w) = [w, y], \; \forall w\in\mathcal A.$$
\end{definition}

\begin{lemma}\label{lemma_W-reduction}
Let $w_{0} = x^{N-k}y^{p}z^{k},$ $w_{i+1}=R_{y}(w_{i})$ where $k, p, N \in\mathbb N_{0}$ and 
$0 \leqslant k \leqslant N,$ Then $w_{N} = N!(-1)^{k}y^{p}$ and $w_{N+i} = 0$ $\forall i\in\mathbb N.$
\end{lemma}
\begin{proof}
Let $k=0.$ Then from lemma~\ref{lemma_AlgebraRelations} it follows:
\begin{align*} 
 w_{1} = R_{y}(w_{0}) &= [w_{0}, y] = [x^{N}y^{p}, y] = Nx^{N-1}y^{p},  \\
 w_{2} = R_{y}(w_{1}) &= N[x^{N-1}y^{p}, y] = N(N-1)x^{N-2}y^p, \\
                      &\ldots \\
 w_{N} = N!y^{p}, \;  & w_{N+1} = R_{y}(w_{N}) = N![y^{p}, y] = 0.
\end{align*}
Case $k=N$ is considered analogously.

Now let $0<k<N.$  Then:
\begin{multline*}
w_{1} =  R_{y}(w_{0}) = [w_{0}, y] = \\
= [x^{N-k}y^{p}z^{k}, y] = x^{N-k}y^{p}z^{k}y - yx^{N-k}y^{p}z^{k} = \\
= x^{N-k}y^{p}(yz^{k}-kz^{k-1} - (x^{N-k}y - (N-k)x^{N-k-1})y^{p}z^{k} = \\
= Nx^{(N-1)-k}y^{p}z^{k} - k(x^{N-k}y^{p}z^{k-1} + x^{N-k-1}y^{p}z^{k})= \\
= Nx^{(N-1)-k}y^{p}z^{k} - kx^{N-k-1}y^{p}(x+z)z^{k-1} - kpx^{N-k-1}y^{p-1}z^{k-1}.
\end{multline*}

First, we have to note that the action of the operator $R_{y}$ onto $w_{i}$ considired as polynom of $x$ 
and $z$ variables decrease its degree by one. Second, acording to lemma~\ref{lemma_AlgebraRelations}
polynoms $(x+z)^{k}$ belong to $Ker R_{y}$ and therefore two last terms in the previous relation may be 
formally regarded as polynoms of degree $N-2$ relatively $x$ and $z.$  So they vanish after $N-1$-times 
action of $R_{y}.$ Summarize, we conclude that $w_{N}$ is a polynom of the only $y$ variable.

In similar way we may obtain:
$$
w_{N-k} = N(N-1)\ldots(k-1)y^{p}z^{k} + \tilde{w}(x, y, z)
$$ 
with $\tilde w(x, y, z)$ denoting terms which vanish after $k$-times action of the operator $R_{y}.$

Continued such process we yield: 
\begin{multline*}
w_{N-k+1} = -N(N-1)\ldots((k-1)ky^{p}z^{k-1} + \tilde{w}(x,y,z), \ldots, \\
w_{N} = (-1)^{k}N!y^{p},\quad w_{N+1} = 0.
\end{multline*}
\end{proof}

\noindent\textsc{Corollary.} \\

 If
 $
 w_{0} = \sum\limits_{p=0}^{M}\left(\alpha_{0 p 0}y^{p} +
         \sum\limits_{k=0}^{1}\alpha_{1-k p k}x^{1-k}y^{p}z^{k}+ \ldots +
         \sum\limits_{k=0}^{N}\alpha_{N-k p k}x^{N-k}y^{p}z^{k}\right)
$ \\
then $w_{N} = N!\sum\limits_{p=0}^{M}\left(\sum\limits_{k=0}^{N}\alpha_{N-k p k}(-1)^{k}\right)y^{p}.$

\begin{lemma}\label{lemma_yP(x+z)SzK-rewrite}
Let $p \in\mathbb N_{0},  N \in\mathbb N$ and $ 0 \leqslant s \leqslant N-1$ then:
\begin{multline*}
w = \sum\limits_{k=0}^{N-s}\gamma_{k}x^{N-s-k}y^{p}(x+z)^{s}z^{k} =\\
 = \delta_{k}  \sum\limits_{k=0}^{N-(s+1)} \Bigl(
x^{N-(s+1)-k}y^{p}(x+z)^{s+1}z^{k} + 
   px^{N-(s+1)-k}y^{p-1}(x+z)^{s}z^{k} -\\
- s(x^{N-(s+1)-k}y^{p+1}(x+z)^{s}z^{k} - x^{N-(s+1)-k}y^{p}(x+z)^{s-1}z^{k})\Bigr) 
\end{multline*}
if $\sum\limits_{k=0}^{N-s}\gamma_{k}(-1)^{k}=0.$
\end{lemma}
\begin{proof}
Let $0 \leqslant k \leqslant N-s-1.$ Using lemma~\ref{lemma_AlgebraRelations} we may write:
\begin{multline*}
x^{N-s-k-1}y^{p}(x+z)^{s+1}z^{k} =
x^{N-s-k-1}y^{p}(x+z)(x+z)^{s}z^{k} =\\
= x^{N-s-k-1}(xy^{p}-py^{p-1})(x+z)^{s}z^{k} + x^{N-s-k-1}y^{p}z(x+z)^{s}z^{k} =\\
= x^{N-s-k}y^{p}(x+z)^{s}z^{k} - px^{N-(s+1)-k}y^{p-1}(x+z)^{s}z^{k} +\\
+ x^{N-(s+1)-k}y^{p}( (x+z)^{s}z + sy(x+z)^{s} + s(x+z)^{s-1} )z^{k}.
\end{multline*}

What can be rewritten as:
\begin{multline*}
x^{N-s-k}y^{p}(x+z)^{s}z^{k} = x^{N-(s+1)-k}y^{p}(x+z)^{s+1}z^{k} - 
x^{N-s-(k+1)}y^{p}(x+z)^{s}z^{k+1} + \\
+ px^{N-(s+1)-k}y^{p-1}(x+z)^{s}z^{k} - x^{N-(s+1)-k}y^{p}(sy(x+z)^{s} + s(x+z)^{s-1})z^{k}.
\end{multline*}

Applying this relation to the first term ($k=0$) of polynom $w$ we obtain:
\begin{multline*}
w = \sum\limits_{k=0}^{N-s}\gamma_{k}x^{N-s-k}y^{p}(x+z)^{s}z^{k} = 
    \gamma_{0}x^{N-(s+1)}y^{p}(x+z)^{s+1} + \\
+ (\gamma_{1}-\gamma_{0})x^{N-s-1}y^{p}(x+z)^{s}z +
 \gamma_{2}x^{N-s-2}y^{p}(x+z)^{s}z^{2} + \ldots + \gamma_{N-s}y^{p}(x+z)^{s}z^{N-s} +\\
+ p\gamma_{0}x^{N-(s+1)}y^{p-1}(x+z)^{s} - \gamma_{0}x^{N-(s+1)}y^{p}(sy(x+z)^{s} + s(x+z)^{s-1}).
\end{multline*}

Then to the monom $x^{N-s-1}y^{p}(x+z)^{s}z:$

\begin{multline*}
w = \gamma_{0}x^{N-(s+1)}y^{p}(x+z)^{s+1} + (\gamma_{1}-\gamma_{0})x^{N-(s+1)-1}y^{p}(x+z)^{s+1}z +\\
+ (\gamma_{2}-\gamma_{1}+\gamma_{0})x^{N-s-2}y^{p}(x+z)^{s}z^{k+1} +\\
+ \gamma_{3}x^{N-s-3)}y^{p}(x+z)^{s}z^{3} + \ldots + \gamma_{N-s}y^{p}(x+z)^{s}z^{N-s} +\\
 + p\left(\gamma_{0}x^{N-(s+1)}y^{p-1}(x+z)^{s} +
(\gamma_{1}-\gamma_{0})x^{N-(s+1)-1}y^{p-1}(x+z)^{s}z\right) - \\
- \bigl(\gamma_{0}x^{N-(s+1)}y^{p}(sy(x+z)^{s} + s(x+z)^{s-1}) + \\
+ (\gamma_{1} - \gamma_{0})x^{N-(s+1)-1}y^{p}(sy(x+z)^{s} + s(x+z)^{s-1})z\bigr).
\end{multline*}

And so on to the monoms $x^{N-s-k}y^{p}(x+z)^{s}z^{k},  k=2,3,\ldots,N-s.$ 
Finally we gain:
\begin{multline*}
 w = \Bigl(\gamma_{0}x^{N-(s+1)}y^{p}(x+z)^{s+1} + 
    (\gamma_{1}-\gamma_{0})x^{N-(s+1)-1}y^{p}(x+z)^{s+1}z + \ldots +\\ 
+ \Bigl[\sum\limits_{k=0}^{N-(s+1)}\gamma_{k}(-1)^{N-(s+1)-k}\Bigr]y^{p}(x+z)^{s+1}z^{N-(s+1)}\Bigr) +\\
+ \Bigl[\sum\limits_{k=0}^{N-s}\gamma_{k}(-1)^{N-s-k}\Bigl]y^{p}(x+z)^{s}z^{N-s} +\\
+ p\Bigl(\gamma_{0}x^{N-(s+1)}y^{p-1}(x+z)^{s} +
(\gamma_{1}-\gamma_{0})x^{N-(s+1)-1}y^{p-1}(x+z)^{s}z + \ldots +\\
+ \Bigl[\sum\limits_{k=0}^{N-(s+1)-k}\gamma_{k}(-1)^{N-(s+1)-k}\Bigr]y^{p-1}(x+z)^{s}z^{N-(s+1)}\Bigr) -\\
- \Bigl(\gamma_{0}x^{N-(s+1)}y^{p}(sy(x+z)^{s} + s(x+z)^{s-1}) + \\
+ (\gamma_{1}-\gamma_{0})x^{N-(s+1)-1}y^{p}(sy(x+z)^{s} + s(x+z)^{s-1})z +\ldots +\\
+ \Bigl[\sum\limits_{k=0}^{N-(s+1)}\gamma_{k}(-1)^{N-(s+1)-k}\Bigr]y^{p}(sy(x+z)^{s} + 
s(x+z)^{s-1})z^{N-(s+1)}\Bigr).
\end{multline*}

Denoting $\gamma_{0}=\delta_{0}, \gamma_{1}-\gamma_{0} = \delta_{1},\ldots, 
\sum\limits_{k=0}^{N-(s+1)}\gamma_{k}(-1)^{N-(s+1)-k}=\delta_{N-(s+1))}$ and taking into account that
$\sum\limits_{k=0}^{N-s}\gamma_{k}(-1)^{N-s-k} = \pm\sum\limits_{k=0}^N\gamma_{k}(-1)^{k} = 0$
we finish the proof.

\end{proof}

\begin{definition}
Define linear operators $L_{x,j} \in End_{\Bfield}(\mathcal A)$ as:\\
$$
 L_{x,j}w = [x, w] + jyw, \forall w \in \mathcal A \mbox{ and } j \in \mathbb N_{0}.
$$
\end{definition}

\begin{lemma}\label{lemma_Ad_xj-action}
 Let $w=y^{i}(x+z)^{k}$ where $i, k\in\mathbb N_{0}.$ Then 
$$
L_{x,j}(w) = iy^{i-1}(x+z)^{k} - ky^{i}(x+z)^{k-1} + (j-k)y^{i+1}y^{i+1}(x+z)^{k}.
$$
\end{lemma}
\begin{proof}

\begin{multline*}
 L_{x,j}(w) =  [x, y^{i}(x+z)^{k}] + jy^{i+1}(x+z)^{k} = 
                   xy^{i}(x+z)^{k} - y^{i}(x+z)^{k}x + jy^{i+1}(x+z)^{k} =\\
= xy^{i}(x+z)^{k} - y^{i}\bigl(x(x+z)^{k} + ky(x+z)^{k} + k(x+z)^{k-1}\bigr) + jy^{i+1}(x+z)^{k} =\\
= xy^{i}(x+z)^{k}- (xy^{i}-iy^{i-1})(x+z)^{k} - ky^{i+1}(x+z)^{k} - ky^{i}(x+z)^{k-1} + jy^{i+1}(x+z)^{k}=\\
= iy^{i-1}(x+z)^{k} - ky^{i}(x+z)^{k-1} + (j-k)y^{i+1}(x+z)^{k}.
\end{multline*}

\end{proof}

\begin{lemma}\label{lemma_Ad_xj-action-result}
Let 
$$
w = \sum_{i=0}^{N}\sum_{k=0}^{M}\alpha_{ij}y^{i}(x+z)^{k} \ne 0, N, M\in\mathbb N_{0}.
$$
where $N = max\{ i \mid \exists \alpha_{ik} \ne 0, k \in\overline{1,M} \},$
$M = max\{ k \mid \exists \alpha_{ik} \ne 0, i \in\overline{1,N} \}.$ 

Then $L_{x,j}(w) \ne 0$ $\forall j\in\mathbb N_{0}$ except the case $N=M=j=0.$
\end{lemma}
\begin{proof}

Consider two cases. The first one is $N=M=0$ i.~e. $w = \alpha_{00} \ne 0.$
then it is clear that $L_{x,j}(w) = 0 \Leftrightarrow j = 0.$

The second one is more general case when $(N,M) \ne (0,0).$ Owing to lemma~\ref{lemma_Ad_xj-action} 
we have:
\begin{multline*}
L_{x,j}(w) = \sum\limits_{i=0}^{N}\sum\limits_{k=0}^{M}\alpha_{ik}L_{x,j}(y^{i}(x+z)^{k}) =\\
 = \sum\limits_{i=0}^{N}\sum\limits_{k=0}^{M}\alpha_{ik}( iy^{i-1}(x+z)^{k} - 
                                                       ky^{i}(x+z)^{k-1} + (j-k)y^{i+1}(x+z)^{k}).
\end{multline*}
From the linear independence of $y^{i}(x+z)^{k}$ it follows that $L_{x,j}(w) \ne 0.$

Indeed, if $w = \sum\limits_{i=0}^{N}\alpha_{iM}y^{i}(x+z)^{M}$ for some $N>0, M \geqslant 0$ then
$$
L_{x,j}(w) = \alpha_{NM}(j-M)y^{N+1}(x+z)^{M} + 
\sum\limits_{i=0}^{N}\sum_{k=M-1}^{M}\tilde{\alpha}_{ik}y^{i}(x+z)^{k} 
$$
and $\alpha_{NM}(j-M) \ne 0$ by $j \ne M.$

Or,
$$
L_{x,M}(w) = \alpha_{NM}(-M) y^{N}(x+z)^{M-1} + 
\sum\limits_{i=0}^{N-1}\sum_{k=M-1}^{M}\tilde{\alpha}_{ik}y^{i}(x+z)^{k} 
$$
and $\alpha_{NM}(-M) \ne 0$ by $ M > 0.$

At $N>0,  j=M=0: L_{x,0}(w) = \sum\limits_{i=0}^{N}\alpha_{i0}iy^{i}.$

Finally denoting $N' = max\{i \mid \alpha_{iM} \ne 0 \}$ we may write:
$$
w = \sum\limits_{i=0}^{N}\sum\limits_{k=0}^{M} \alpha_{iM}y^{i}(x+z)^{k} =
    \sum\limits_{i=0}^{N'}\alpha_{iM} y^{i}(x+z)^{M} + 
    \sum\limits_{i=0}^{N}\sum\limits_{k=0}^{M-1} \alpha_{iM}y^{i}(x+z)^{k}
$$
and prove as previously $L_{x,j}(w) \ne 0.$

\end{proof}

\begin{proposition}\label{prop_simplicity_A}
 $\mathcal A$ is a simple algebra.
\end{proposition}
\begin{proof}

Let $0 \ne w\in I \trianglelefteq\mathcal A$ and 
$$
w = \sum\limits_{p=0}^{M}\left(\alpha_{0 p 0}y^{p} + 
    \sum\limits_{k=0}^{1}\alpha_{1-k p k}x^{1-k}y^{p}z^{k} + \ldots +
    \sum\limits_{k=0}^{N}\alpha_{N-k p k}x^{N-k}y^{p}z^{k}\right).
$$
Consider two cases.

1) Suppose that  $w = \sum\limits_{p=0}^{M}\alpha_{0 p 0}y^{p}, \alpha_{0 M 0} \ne 0.$
Then $(L_{x,0})^{M}(w) = M!\alpha_{0 M 0} = \tilde{w} \ne 0 $ and $\tilde{w} \in I.$

2) Now show that general case may be reduced to the first special one. Applying recurrently
lemma~\ref{lemma_yP(x+z)SzK-rewrite} if it is need we may always represent $w$ for some 
$K: 0 \leqslant K \leqslant N$ as:
$$
w = \sum\limits_{p=0}^{M}\sum\limits_{s=0}^{N-K} \left(
\sum\limits_{k=0}^{K}\beta_{k}^{ps}x^{K-k}y^{p}(x+z)^{s}z^{k} +  
\sum\limits_{k=0}^{K-1}\gamma_{k}^{ps}x^{K-1-k}y^{p}(x+z)^{s}z^{k} + \ldots \right) 
$$
where $\exists p \in\overline{0,M}: \sum\limits_{k=0}^{K}\beta_{k}^{ps}(-1)^{k} \ne 0.$

Exactly as it was made in the proof of lemma~\ref{lemma_W-reduction} it can be shown that:
$$
\bar{w} = (R_{y})^{K}(w) = K!\sum\limits_{p=0}^{M}\sum\limits_{s=0}^{N-K}
 \left(\sum\limits_{k=0}^{K}\beta_{k}^{ps}(-1)^{k}\right) y^{p}(x+z)^{s} \in I.
$$

Note that $\bar{w} \ne 0$ due to the linear independence of $y^{p}(x+z)^{s}$ (not all 
$\sum\limits_{k=0}^{L}\beta_{k}^{ps}(-1)^{k}$ equals 0).

Rewrite $\bar w$ as $\bar w = \sum\limits_{i=0}^{\bar{M}}\sum\limits_{k=0}^{\bar{N}}
\bar{\alpha}_{ik}y^{i}(x+z)^{k},0 \leqslant \bar{M}\leqslant {M}, 0 \leqslant \bar{N} \leqslant N-K .$ 
(without loss of generality (w.~l.~o.~g) one may be assumed that $\bar{\alpha}_{\bar{M}\bar{N}}\ne 0.$)

Let $\tilde{w} = (L_{x,\bar {N}})^{\bar{M}+1}(\bar{w})$ then owing to 
lemmas~\ref{lemma_Ad_xj-action} and~\ref{lemma_Ad_xj-action-result} we have 
$\deg_{x+z}\tilde{w}\leqslant \bar{N}-1$ and $\tilde{w}\ne 0.$ Applying operators
$(L_{x,i_{k}})^{m_{k}}, i_{k} \ne 0$ after some steps we get:
$w' = (L_{x,i_{s}})^{m_{s}}\ldots (L_{x,i_{1}})^{m_{1}}(\tilde{w}) \ne 0, w'\in I$ and 
$\deg_{x+z} w' = 0$ i.~e. $w'= \sum\limits_{k=0}^{M'}\alpha_{k}y^{k} \ne 0, w' \in I.$

Summing two cases we gain: if $0 \ne I \trianglelefteq\mathcal A$ then 
$1_{\mathcal A} \in I \Rightarrow I = \mathcal A$ i.~e. \\
$\mathcal A$ is a simple algebra.

\end{proof}

\begin{definition}
Fix $n \in \mathbb N.$ Define an algebra $\mathcal A_{n}$ over field $\Bfield$ as algebra 
for whose generators $\{ x_{1},x_{2},\ldots,x_{n},y,z_{1},z_{2},\ldots,z_{n} \}$  the following
relations hold:
\begin{align*}
&[x_{1}, y]     = [y, z_{1}] = [x_{i}, z_{i}] = 1, \;i = 2,\ldots, n; \\
&[x_{1}, z_{1}] = -yz_1-xy_1;  \\
&[x_{i}, x_{j}] = [x_{i}, z_{j}] = [z_{i}, z_{j}] = 0, \; i \ne j, i, \; j = 1,\ldots, n; \\
&[x_{i}, y]     = [y, z_{i}] = 0, \; i=2, \ldots, n.
\end{align*}
\end{definition}

Show that $\mathcal A_{n}$ is a simple algebra with $\GKdim(\mathcal A_{n}) = 2n+1.$
We shall prove this for the case $n=2$ since a general case is considered analogously.

As well as the case $n=1$ it may be shown that monoms $x_{1}^{i}x_{2}^{j}y^{p}z_{1}^{k}z_{2}^{l}$
form a basis of algebra $\mathcal A_{2}$ over $\Bfield$ and particulary $\GKdim(\mathcal A_{2})=5.$ 

First define linear operators $L_{x_{2}}$ and $R_{z_{2}}\in End_{\Bfield}(\mathcal A_{2})$ 
such as:
$$
 L_{x_{2}}(w)=[x_{2}, w], R_{z_{2}}(w)=[w, z_{2}], \forall w\in \mathcal A_{2}.
$$

Let now
$$
w_{0} = \sum\limits_{i=0}^{N_{1}}\sum\limits_{j=0}^{N_{2}}\sum\limits_{p=0}^{N_{3}}
\sum\limits_{k=0}^{N_{4}}\sum\limits_{l=0}^{N_{5}}\alpha_{ijpkl}x_{1}^{i}x_{2}^{j}y^{p}z_{1}^{k}z_{2}^{l}
\ne 0
$$ and $w_{0} \in I \trianglelefteq\mathcal A.$ W.l.o.g. one may assume that $\deg_{x_{2}}w = N_{2}, \deg_{z_{2}}w = N_{5}.$ Then it might be obtained 
(exactly as for lemma~\ref{lemma_AlgebraRelations}) that:
\begin{gather*}
w_{1} = (L_{x_{2}})^{N_{5}}(w_{0}) = 
N_{5}!\sum\limits_{i=0}^{N_{1}}\sum\limits_{j=0}^{N_{2}}\sum\limits_{p=0}^{N_{3}}\sum\limits_{k=0}^{N_{4}}\alpha_{ijpkN_{5}}x_{1}^{i}x_{2}^{j}y^{p}z_{1}^{k} \ne 0, w_{1} \in I. \\
w_{2}= (R_{z_{2}})^{N_{2}}(w_{1}) = N_{2}!N_{5}!\sum\limits_{i=0}^{N_{1}}\sum\limits_{p=0}^{N_{3}}\sum\limits_{k=0}^{N_{4}}\alpha_{iN_{2}pkN_{5}}x_{1}^{i}y^{p}z_{1}^{k} \ne 0, w_2\in I.
\end{gather*}

Verbatim repeating the proof of proposition~\ref{prop_simplicity_A} we yield: 
$1_{\mathcal A_{2}} \in I \Rightarrow \mathcal A_{2}$ is a simple algebra.

\end{document}